\documentclass[12pt]{article}
\usepackage{soul,color}
\newcommand{\tr}{^{\prime}}
\def\b#1{\mbox{\boldmath $#1$}}    
\def\bl#1{\mbox{\footnotesize \boldmath {$#1$}}} 
\def\cg#1{\mbox{${\cal #1}$}}

\newcommand{\ot}{\mbox{$\:\otimes \:$}}

\begin{document}
\title{\vspace*{-1cm}Bayesian inference 
through encompassing priors and importance sampling\\ for a class of
marginal models for categorical data}
\author{Francesco Bartolucci\footnote{Dipartimento di Economia, Finanza e
Statistica, Universit\`{a} di Perugia, Via A. Pascoli 20, 06123
Perugia, Italy, {\em email}: bart@stat.unipg.it}\hspace{1mm}, Luisa
Scaccia\footnote{Dipartimento di Istituzioni Economiche e
Finanziarie, Universit\`a di Macerata, Via Crescimbeni 20, 62100
Macerata, Italy, {\em email}: scaccia@unimc.it}\hspace{1mm} and
Alessio Farcomeni\footnote{Dipartimento di Sanit\'{a}
Pubblica e Malattie Infettive,
Sapienza - Universit\`a di Roma, Piazzale Aldo Moro, 5,
00186 Roma, Italy, {\em email}: alessio.farcomeni@uniroma1.it} }
\maketitle

\vspace*{-0.5cm}

\begin{abstract}
\noindent We develop a Bayesian approach for selecting the model
which is the most supported by the data within a class of
marginal models for categorical variables formulated through
equality and/or inequality constraints on generalised logits (local,
global, continuation or reverse continuation), generalised log-odds
ratios and similar higher-order interactions. For each constrained
model, the prior distribution of the model parameters is formulated
following the encompassing prior approach. Then, model selection is
performed by using Bayes factors which are estimated by an
importance sampling method. The approach is illustrated through
three applications involving some datasets, which also include
explanatory variables. In connection with one of these examples, a
sensitivity analysis to the prior specification is also considered.

\vspace{0.3cm} \noindent {\bf Keywords}: Bayes factor; Encompassing priors;
Generalised logits; Inequality constraints; Marginal Likelihood; Positive association.

\end{abstract}\newpage
\section{Introduction}
Though log-linear models are frequently used for the analysis of
contingency tables, they do not allow to express, and consequently
test, several hypotheses that are usually of interest, mainly
because lower order interactions do not refer to the marginal
distributions to which they seem to refer. This motivated McCullagh
\& Nelder (1989, Section 6.5) to introduce a class of models in
which the joint distribution of a set of categorical variables is
parametrised through the highest log-linear interactions within each
possible marginal distribution. Several other models have been
proposed following the original idea of McCullagh \& Nelder (see
Glonek \& McCullagh, 1995, Glonek, 1996, Colombi \& Forcina, 2001,
and Bergsma \& Rudas, 2002, and Bartolucci {\em et al.}, 2007).

In this paper, we deal with a flexible class of models in which: (i)
the parameters of the saturated model are given by {\em generalised}
logits, in the sense of Douglas {\em et al.} (1991), for each
univariate marginal distribution, generalised log-odds ratios for
each bivariate marginal distribution and similar interactions for
each higher-order marginal distribution; (ii) any constrained model
may be formulated through linear equality and inequality constraints
on such parameters. In this way we may express several hypotheses
which are of special interest in presence of ordinal variables (see
Bartolucci {\em et al.}, 2001, and Colombi \& Forcina, 2001), as for
instance, that: (i) the marginal distribution of one variable is
stochastically larger than that of another variable, provided that
these have the same categories; (ii) a certain type of positive
association between a pair of variables holds; (iii) the marginal
distribution of one variable is stochastically increasing with
respect to the level of an explanatory variable.

For the above class of models we develop Bayesian inference and, in
particular, a model selection strategy based on the Bayes factor
(see Jeffreys, 1935, 1961, and Kass \& Raftery, 1995), which is
defined as the ratio between the marginal likelihoods of two
competing models. For the marginal models considered in this paper,
the use of the Bayes factor allows us to easily compare two models
parametrised through different types of logit, which would be
otherwise cumbersome using a likelihood ratio test. Moreover, since
the Bayes factor is computed as the ratio between two marginal
likelihoods, the presence of the nuisance parameters does not affect
the inferential results, as it typically happens when using a
likelihood ratio test (for a discussion on this point see Dardanoni
\& Forcina, 1998, and Bartolucci {\em et al.}, 2001). On the other
hand, the proposed approach requires the specification of a prior
distribution on the parameters that is not required within the
likelihood ratio approach.

While the decision theoretic approach leads us to select the model
with largest marginal likelihood, we can also use the Bayes factor as
a measure of evidence. In order to assess this evidence we refer to
the Jeffreys (1961) scale, which gives the following guideline: a log
Bayes factor below 0.5 indicates {\it poor} evidence, between 0.5 and
1 {\it substantial}, between 1 and 2 {\it strong} and {\it decisive}
evidence is provided by a Bayes factor larger than 2. 

Bayesian methods for the analysis of categorical data has been dealt
with by several authors. For instance, Albert (1996, 1997) used the
Bayes factor to test hypotheses such as independence,
quasi-independence, symmetry or constant association in two-way and
three-way contingency tables. Dellaportas \& Forster (1999) proposed
a general framework for selecting a log-linear model through the
Reversible Jump algorithm of Green (1995) under a multivariate
Normal prior distribution on the parameters. In practice, both
Albert (1996, 1997) and Dellaportas \& Forster (1999) dealt with
log-linear models obtained by imposing some linear equality
constraints on the parameters of the saturated model, for example
that a subset of the parameters is equal to zero. Klugkist and
Hoijtink (2007), Hoijtink {\em et al.} (2008), Klugkist {\em et al.}
(2005a, 2005b, 2010) and Wetzels {\em et al.} (2010) used, instead,
the Bayes factor to compare competing models expressed through
linear inequality and {\em about equality} constraints on the saturated
model. Under their \textit{encompassing prior} approach, the Bayes
factor between a constrained model and the encompassing model
reduces to the ratio of the probability that the constraints hold
under the encompassing posterior distribution and the probability that
they hold under the encompassing prior distribution. By encompassing
model we mean a model, the parameter space of which includes that of
every other model under consideration. Therefore, once the prior
distribution has been specified on the encompassing model parameters
(encompassing prior), it is automatically specified for each
submodel.

The selection strategy we adopt for the class of models considered
in this paper is related to the approach of Klugkist {\it et al.}
(2010). We exploit their encompassing prior approach, which leads to a
logically coherent assessment of prior and posterior model
probabilities and parameters distributions, as well as an easy
estimation of the Bayes factors. However, our work differs from that
of Klugkist {\it et al.} (2010) mainly in three respects: (i) we
consider a more general class of models for categorical data; (ii)
we propose an importance sampling method to improve the efficiency
of the Bayes factor estimates for models with very small prior and,
possibly, posterior probabilities; (iii) we introduce an iterative
algorithm to estimate the Bayes factor for models specified through
about equality constraints, which does not require to sample
from a constrained model parameter space.

The paper is organised as follows. In Section 2 we describe the
class of models of interest. Then in Section 3 we review the
encompassing prior approach and we deal with Bayesian model
selection. Finally, in Section 4 we illustrate the proposed approach
through three applications involving some datasets of interest in
the categorical data analysis literature.
\section{Marginal models for categorical variables}
In this section, we introduce the marginal models developed by
McCullagh \& Nelder (1989) and illustrate the parameterisation based
on generalised logits and log-odd ratios. Then, we show how
hypotheses of interest may be expressed through linear equality and
inequality constraints imposed on the parameters of the saturated
model.
\subsection{Parametrisation}
Let $\b A=(A_1,\ldots,A_q)$ be a vector of $q$ categorical variables
and $\{1,\ldots,m_i\}$ be the support of $A_i$. Also let
$r=\prod_im_i$ be the number of possible configurations of $\b A$
and $\b \pi$ be the $r$-dimensional column vector of the joint
probabilities $p(\b A=\b a)$ arranged in lexicographical order. In
the following, we describe a saturated parameterisation of such a
vector based on marginal logits, marginal log-odds ratios and
similar higher-order interactions.

Marginal logits may be of type {\em local} ($l$), {\em global}
($g$), {\em continuation} ($c$) or {\em reverse continuation}
($r$). For the $i$-th variable these are defined as follows for
$a_i = 1,\ldots,m_i-1$:
\begin{itemize}
\item {\em local}:\hspace{3.55cm} ${\displaystyle \eta_i(a_i;l) =
\log\frac{p(A_i=a_i+1)}{p(A_i=a_i)}}$;
\item {\em global}:\hspace{3.35cm} ${\displaystyle \eta_i(a_i;g) =
\log\frac{p(A_i\geq a_i+1)}{p(A_i\leq a_i)}}$;
\item {\em continuation}:\hspace{2.1cm} ${\displaystyle
\eta_i(a_i;c) = \log\frac{p(A_i\geq a_i+1)}{p(A_i= a_i)}}$;
\item {\em reverse continuation}: \hspace{0.5cm} ${\displaystyle
\eta_i(a_i;r) = \log\frac{p(A_i=a_i+1)}{p(A_i\leq a_i)}}$.
\end{itemize}
Local logits are used when it is of interest to compare the marginal
probability of each category with that of the previous category.
Logits of type global and continuation are specially tailored to
ordinal variables. In particular, logits of type global are more
appropriate when the variable may be seen as a discretised version
of an underlying continuum, whereas logits of type continuation are
more appropriate when categories correspond to levels of achievement
that may be entered only if the previous level has also been
achieved, as in education. Finally, using logits of type reverse
continuation is the same as arranging categories in reverse order
and using logits of type continuation.

Marginal log-odds ratios are defined as contrasts between
conditional logits. For two variables, $A_i$ and $A_j$, the most
well-known log-odds ratios are shown in the following, where
$a_i=1,\ldots,m_i-1$ and $a_j=1,\ldots,m_j-1$:
\begin{itemize}
\item {\em Local}: when logits of type $l$ are
used for both $A_i$ and $A_j$
\begin{eqnarray*}
\eta_{ij}(a_i,a_j;l,l) &=&
\eta_j(a_j;l|A_i=a_i+1)-\eta_j(a_j;l|A_i=a_i) =\\
&=&
\log\frac{p(A_i=a_i,A_j=a_j)p(A_i=a_i+1,A_j=a_j+1)}{p(A_i=a_i+1,A_j=a_j)p(A_i=a_i,A_j=a_j+1)};
\end{eqnarray*}
\item {\em Local-Global}: when logits of type $l$ are used for
$A_i$ and of type $g$ for $A_j$
\begin{eqnarray*}
\eta_{ij}(a_i,a_j;l,g) &=&
\eta_j(a_j;g|A_i=a_i+1)-\eta_j(a_j;g|A_i=a_i) =\\
&=& \log\frac{p(A_i=a_i,A_j\leq a_j)p(A_i=a_i+1,A_j\geq
a_j+1)}{p(A_i=a_i+1,A_j\leq a_j)p(A_i=a_i,A_j\geq a_j+1)};
\end{eqnarray*}
\item {\em Global}: when logits of type $g$ are used for both
$A_i$ and $A_j$
\begin{eqnarray*}
\eta_{ij}(a_i,a_j;g,g) &=&
\eta_j(a_j;g|A_i\geq a_i+1)-\eta_j(a_j;g|A_i\leq a_i) =\\
&=& \log\frac{p(A_i\leq a_i,A_j\leq a_j)p(A_i\geq a_i+1,A_j\geq
a_j+1)}{p(A_i\geq a_i+1,A_j\leq a_j)p(A_i\leq a_i,A_j\geq a_j+1)}.
\end{eqnarray*}
\end{itemize}
Similarly, three-way interactions are defined as contrasts between
conditional log-odds ratios and so on for higher order
interactions.

Now let $\b z$ be an $r$-dimensional vector of zeros and ones, let
$\b\eta_{\bl z}$ be a column vector containing all the marginal
interactions between the variables $A_i$ such that $z_i = 1$ and let
$\b\eta$ be the vector obtained by stacking, in lexicographical
order, the vectors $\b\eta_{\bl z}$ one below the other for any $\b
z\neq \b 0$. Following Colombi \& Forcina (2001), such a vector,
which provides the saturated parametrisation of $\b\pi$ at issue, may
be simply obtained as
\begin{equation}
\b\eta = \b C\log(\b M\b\pi), \label{marg}
\end{equation}
where $\b C$ and $\b M$ are appropriate matrices, whose construction
is described in Appendix A. Note, however, that to
invert equation (\ref{marg}), and so obtain $\b\pi$ in terms of
$\b\eta$, we must rely on a Newton-Raphson algorithm as the one
described in Glonek \& McCullagh (1995) and Colombi \& Forcina
(2001).
\subsection{Constrained models}\label{models}
A variety of constrained models may be formulated by posing
linear equality and inequality constraints of the form
\begin{equation}
\b E\b \eta = \b 0,\quad \b U\b \eta\geq \b 0,\label{con}
\end{equation}
on the saturated parameter vector. Here, and throughout the paper,
equality constraints are substituted by
about equality constraints of type $|\b E\b \eta| \leq \b \epsilon$,
for a vector $\b\epsilon>\b 0$ having suitably small elements; this
choice is motivated in Section \ref{compbf}.

Consider first the case of only two variables, $A_1$ and $A_2$. The
most interesting hypotheses are usually on the association between
these variables. Let $c=m_1+m_2-2$ be the number of the marginal
logits and $d=(m_1-1)(m_2-1)$ be that of the log-odds ratios. By
requiring that all the log-odds ratios are non-negative, namely by
letting
\begin{equation}
\b U = \pmatrix{\b O_{d,c} & \b I_d},\label{ind}
\end{equation}
where $\b O_{d,c}$ is a matrix of $d\times c$ zeros and $\b I_d$
denotes an identity matrix of dimension $d$, we express the
hypothesis of positive association between $A_1$ and $A_2$.
Obviously, the type of association depends on the type of logit that
is used for the two variables. For instance, with local logits for
both variables we are formulating the hypothesis of {\em Total
Positivity of Order 2} (TP$_2$; see Karlin, 1968), whereas with
global logits for both variables we are formulating the hypothesis
of {\em Positive Quadrant Dependence} (PQD; see Lehmann, 1966). Note
that there is a hierarchy among these notions of positive
association in the sense that, for instance, TP$_2$ implies that all
the continuation log-odds ratios are non-negative which, in turn,
implies PQD (see Douglas {\em et al.} 1991, for details). Also note
that, regardless of the type of log-odds ratio, independence between
$A_1$ and $A_2$ may be expressed through the constraint that $\b E$,
rather than $\b U$, is equal to the matrix in (\ref{ind}). A less
stringent constraint than that of independence is the constraint of
uniform association, namely that all the $d$ log-odds ratios are
equal to each other (Plackett, 1965). This type of constraint is
formulated by letting
\[
\b E = \pmatrix{\b O_{d-1,c} & \b D_{d-1}},
\]
where, in general, $\b D_h = \pmatrix{\b 0_{h-1} & \b I_{h-1}} -
\pmatrix{\b I_{h-1} & \b 0_{h-1}}$, with $\b 0_h = \b O_{h,1}$, is
a matrix that produces first differences.

When $A_1$ and $A_2$ have the same categories, the number of which is indicated by $m$,
constraints on the univariate marginal distributions may also be of
interest. For instance, we may formulate the constraint of marginal
homogeneity by letting
\[
\b E = \pmatrix{-\b I_{m-1} & \b I_{m-1} & \b O_{m-1,d}}.
\]
When global logits are used for both variables, the constraint that
$A_2$ is stochastically greater than $A_1$ may be imposed by letting
$\b U$, rather than $\b E$, equal to the matrix above. When there
are more than two variables, similar constraints may also involve
interactions of order higher than two. These are typically of
interest in the presence of longitudinal data. Some useful insights
on how formulating a marginal model in these situations are provided
by Glonek \& McCullagh (1995, Section 5).

The approach outlined above may be easily extended when dealing with
one or more explanatory variables, which are collected in 
vector $\b B$. Every possible configuration of these variables, say
$\b b$, defines a stratum conditionally on which we have a vector of
joint probabilities of the response variables, $\b \pi(\b b)$, and a
vector of marginal parameters, $\b\eta(\b b)$, defined as in
(\ref{marg}). Obviously, in this setting we may impose the same
constraints illustrated above within each stratum and also
constraints involving the parameters of different strata. These
constraints are still of the type $\b E\b\eta=\b 0$, $\b U\b\eta\geq
\b 0$, but in this case by $\b\eta$ we mean the vector obtained by
stacking one below the other the vectors $\b\eta(\b b)$ for every $\b
b$. For instance, when we only have two response
variables, we may express the constraint of conditional
independence between $A_1$ and $A_2$, given $\b B$, as
\[
\b E = \b I_s\ot\pmatrix{\b O_{d,c} & \b I_d},
\]
where $s$ is the number of strata, that is the number of different configurations
of $\b B$, and $\ot$ denotes the Kronecker product. We may also
formulate the hypothesis that $A_1$ and $A_2$ have the same degree
of association for each stratum, by letting
\[
\b E=\b D_s\ot\pmatrix{\b O_{d,c} & \b I_d},
\]
or that the explanatory variables do not affect the marginal
distributions of $A_1$ and $A_2$, by letting
\[
\b E=\b D_s\ot\pmatrix{\b I_c & \b O_{c,d}}.
\]
Finally, if we have only one explanatory variable, $B$, and this
is ordinal, we can express the constraint that the marginal
distributions of $A_1$ and $A_2$ increase with the level of $B$ by
letting $\b U$, rather than $\b E$, equal to the matrix above.
\section{Bayesian estimation and model selection}
In this section we show how to make inference on the models
presented in the previous section. In particular, Section
\ref{priors} is devoted to the issue of choosing appropriate priors
for the class of models at hand, whereas Section \ref{modsel} is
focused on assessing the plausibility of the different models, given
an observed contingency table.

In the following, when the data are not stratified, we denote the
frequency corresponding to the configuration $\b a$ of such a table
by $y_{\bl a}$ and by $\b y$ the vector with elements $y_{\bl a}$
arranged as in $\b \pi$. When the data are stratified according to
one or more explanatory categorical variables, we have a vector of
frequencies $\b y(\b b)$ for every configuration of such variables
and, consequently, $\b y$ denotes the vector obtained by stacking
one below the other the vectors $\b y(\b b)$ for every ${\b b}$.
\subsection{Prior distributions}\label{priors}
In a Bayesian framework, it is natural to include equality and
inequality constraints imposed on the model parameters as prior
knowledge. Since all the constrained models presented in Section
\ref{models} are nested in an unconstrained or encompassing model,
we use the concept of encompassing prior (Klugkist {\it et al.},
2005a, 2005b; Klugkist \& Hoijtink, 2007). Therefore, we specify the
prior distribution only for an encompassing model and then we derive
the prior distributions for the other models by restricting the
parameter space according to the constraints of interest. This
approach has the very nice interpretation that the resulting Bayes
factor for model selection (see Section \ref{modsel}) coincides with
the ratio between the proportions of the parameter space that are in
agreement with the constrained model, under the posterior and the
prior distributions of the encompassing model. This approach also
has the advantage that only one single prior distribution needs to be specified.
Moreover, the method can be seen as a generalisation of the
Savage-Dickey density ratio, which overcomes the Borel-Kolmogorov
paradox (Dawid \& Lauritzen, 2001). See Wetzels {\it et al.} (2010)
for a detailed discussion on this point.

In the present framework, it is natural to choose the saturated
model based on the parameter vector $\b\pi$ as the encompassing
model. Under this model, the parameter space is the simplex of
dimension $r$ and the frequency vector $\b y$ has multinomial
distribution with parameters $n$ and $\b\pi$. The choice of the
$\b\pi$ parameterisation, rather than the parameterisation based on
the vector of marginal parameters $\b\eta$, is motivated by the fact
that it also makes straightforward the comparison between different
types of logit.

Let $M_1$ indicate the saturated (encompassing) model and let
$p(\b\pi|M_1)$ denote the encompassing prior distribution. The prior
distribution of each constrained model $M_k$, for $k = 2, \ldots,
K$, follows directly from this prior as
\begin{equation}\label{priorcon} p(\b \pi |M_k) =
\frac{p(\b \pi |M_1) \delta_k(\b\pi)}
{\int p(\b \pi |M_1) \delta_k(\b\pi) \mbox{d}\b \pi} = c_k p(\b \pi |M_1)
\delta_k(\b\pi),\end{equation}
where the integral is on the simplex of dimension $r$. Moreover,
$\delta_k(\b\pi)$ is the indicator function equal to $1$ if $\b\pi$ is
in accordance with the constraints defining model $M_k$ and to $0$ otherwise
and $c_k$ is the inverse of the proportion of the parameter space that, under the
encompassing prior, is in agreement with these constraints.
Obviously, the constrained prior in (\ref{priorcon}) is not defined for a model
with equality constraints, but it is defined for a model with about equality 
constraints.

Under the encompassing prior approach, also the posterior
distribution of the parameters for each constrained model
immediately follows from the posterior under the encompassing model.
In particular, we have
\begin{equation}\label{postcon} p(\b \pi |\b y, M_k) = \frac{p(\b \pi |\b y, M_1)
\delta_k(\b\pi)}{\int p(\b\pi|\b y,M_1) \delta_k(\b\pi) \mbox{d}\b \pi} = d_k p(\b \pi |\b y, M_1) \delta_k(\b\pi),\end{equation}
where, now, $d_k$ is the inverse of the proportion of the
parameter space that, under the encompassing posterior, is in agreement with the constraints of
model $M_k$.

Coming to the issue of choosing a distributional shape for the encompassing prior 
$p(\b \pi |M_1)$, the default prior for $\b\pi$ has been
acknowledged to be the one in which $\b\pi$ has a uniform
distribution on the simplex of dimension $r$ or, equivalently, $p(\b
\pi|M_1)\sim D(\b 1_r)$, where $D(\cdot)$ denotes the Dirichlet
distribution and $\b 1_r$ is a column vector of $r$ ones.
See, for instance, Tuyl {\it et al.} (2009) for a detailed discussion on
this choice.

The posterior for the saturated parameterisation $\b\pi$, with the
default prior choice, is readily derived and it is of type $D(\b
1_r+\b y)$. Therefore, samples can be drawn independently from the
prior and posterior distributions for the saturated model and the
corresponding normalising constants are available in closed form.
\subsection{Model Selection}\label{modsel}
Let $\cg M = \{M_1, \ldots, M_K\}$ denote the set of models of
interest. As already noted, each of these models is defined by a
certain type of logit for every response variable and by constraints
of type (\ref{con}) on the vector of marginal parameters, with the
exception of model $M_1$ which is the saturated model. Then, $M_2,
\ldots, M_K$ are all nested in $M_1$, but not necessarily nested
in one another.

For model selection, we make use of the {\em Bayes factor}, which is
the ratio of the marginal likelihoods of two competing models. Thus,
the Bayes factor for model $M_k$ versus the encompassing model is
defined as:
$$B_{k1} = \frac{p(\b y|M_k)}{p(\b y|M_1)}=\frac{\int p(\b y|\b \pi,M_k)
p(\b \pi|M_k)\mbox{d}\b \pi}
{\int p(\b y|\b \pi,M_1)p(\b \pi|M_1)\mbox{d}\b \pi},$$ where $p(\b
y|\b \pi,M_k)$ and $p(\b y|M_k)$ denote, respectively, the
likelihood of the data and the marginal likelihood for model $M_k$.
The Bayes factor measures the evidence that the data provide for one
model versus the other and correponds to the fold change from
prior model odds to posterior model odds. In this paper we always use
a 0-1 loss. 
Obviously, the larger is $B_{k1}$, the greater is the evidence
provided by the data in favour of $M_k$ with respect to $M_1$ (see
Kass \& Raftery, 1995). So, when $B_{k1}$ is larger than 1, or
equivalently $\log(B_{k1})>0$, model $M_k$ has to be preferred
to model $M_1$. To compare more than two models, or equivalently to choose the
best model in $\cg M$ when $K>2$, a convenient possibility is to
single out $M_1$ as the reference model
and then compute the Bayes factor between every other model and the
unconstrained one, that is $B_{k1}$, for $k=2,\ldots,K$. The model
to be preferred is that with the largest Bayes factor, provided that
it is larger than 1; otherwise the best model is the saturated
model. Obviously, the Bayes factor for comparing every pair of
models $M_k$ and $M_l$, not necessarily nested, is straightforwardly
computed as $B_{kl}=B_{k1}/B_{l1}$.

It is important to note that the Bayes factor, as model selection
tool, combines goodness of fit with a correction for model complexity.
\subsubsection{Computational issues in estimating the Bayes factor}\label{compbf}
Direct computation of the Bayes factor is almost always infeasible,
and this also happens for the class of models dealt with here.
Several methods have been proposed to estimate the Bayes factor
numerically, but the estimation is generally cumbersome from the
computational point of view. 

The encompassing prior approach renders
a nice interpretation of the Bayes factor for a constrained model
$M_k$ with the encompassing model $M_1$, which virtually eliminates
the computational complications inherent in Bayes factor estimation.
In fact, as demonstrated in Klugkist {\it et al.} (2005a), the Bayes factor
for a constrained versus the encompassing model reduces to the ratio
of the proportions of the parameter space that are in agreement with
the constrained model under the posterior distribution and prior distribution of
the encompassing model. Thus, the Bayes factor for a
constrained model $M_k$ with respect to the encompassing model $M_1$ is
\begin{equation}\label{bf3}
B_{k1}=\frac{c_k}{d_k}.
\end{equation}

In the light of (\ref{bf3}), estimating the Bayes factor is particularly simple. 
The encompassing prior is sampled and
$c_k$ is estimated by $\hat c_k$, which is the inverse of the proportion of the 
sample 
that is in agreement with the constraints defining model $M_k$. Similarly, 
sampling
from the encompassing posterior allows us to estimate $d_k$ as $\hat d_k$, which 
is the inverse
of the proportion of the sample that is in agreement with the constraints of model
$M_k$. In this way, using just one sample from the encompassing prior and 
another from the encompassing posterior, the estimate
$$\hat B_{k1}=\frac{\hat c_k}{\hat d_k}$$ 
can be computed for each constrained model $M_k$, $k=2,\ldots,K$.

Notice that, in our setting, the choice of the Dirichlet default
prior for $\b \pi$ allows to sample independently under both the
encompassing prior and posterior distributions, leading to further
simplifications in estimating the Bayes factor. Moreover, in some
cases, $c_k$ can be computed exactly, without the need of sampling
from the encompassing prior. However, there are two issues that we
must deal with when estimating the Bayes factor.

First of all, a rare event problem can arise. Consider for
instance our example in Section \ref{firstexample}: we have a six by
six contingency table with 35 free parameters under the
unconstrained model. The hypothesis of positive association is formulated
by requiring the positivity of the
25 log-odds ratios. When using
logits of type $l$ for both variables, the constant $c_k^{-1}$
for a positive association model can be
calculated exactly as $0.5^{25} = 2.9802\times 10^{-8}$. In this
case, sampling from the encompassing prior is not required, but
such a small values of $c_k^{-1}$
can be common to other models. For these models, even if we drew 
millions of values from the encompassing prior, we would expect to see no
values satisfying the constraint. This would lead an estimate of
the Bayes factor equal to $\infty$ or to $\infty/\infty$, in case
the same problem also arises when sampling from the encompassing
posterior. In general, even if a finite estimate of the Bayes factor
can be achieved, its variance would be huge for those constrained
models characterised by a very small proportion of the parameter space in
agreement with the constraints under the encompassing prior and,
possibly, posterior distribution.

The problem described above is that of {\it rare event simulation}
 (e.g., Bucklew, 2004), which is often overcome 
through \textit{importance sampling}.
Suppose we want to estimate $1/c_k$. From (\ref{priorcon}) it
immediately follows that $1/c_k= E_p(\delta_k(\b\pi))$, where the
expected value is calculated with respect to the encompassing prior
$p(\b\pi|M_1)$. Now, letting $g(\b\pi)$ be any other density such
that $p(\b\pi|M_1)=0$ whenever $g(\b\pi)=0$, we can re-write
\begin{equation}\label{impsamp}
E_p(\delta_k(\b\pi)) = \int p(\b \pi |M_1) \delta_k(\b\pi) \mbox{d}\b \pi =
\int \left[\frac{p(\b \pi |M_1)}{g(\b\pi)} \delta_k(\b\pi)\right]g(\b\pi)
\mbox{d}\b \pi = E_g\left[\delta_k(\b\pi)\frac{p(\b \pi |M_1)}{g(\b\pi)}\right],
\end{equation}
where the last expected value is now calculated under the
\textit{importance} density $g(\b\pi)$. Then, an importance sampling
estimate of $1/c_k$ can be obtained by sampling $\b\pi$ from
an appropriate importance density $g(\b\pi)$ and estimating the last
expected value in (\ref{impsamp}) through the sample mean. If
required, an estimate of $1/d_k$ can be obtained in a similar way, after
choosing an appropriate sampling density.

In this paper, we propose an automatic way to obtain an adequate
importance sampling density. Suppose we want to estimate $1/d_k$ for
a certain model $M_k$; then the proposed method is based on the
following steps:
\begin{enumerate}
\item[(i)] compute the maximum likelihood estimate of the vector
$\b\eta$ under the constraints imposed by model $M_k$ (see Colombi
\& Forcina, 2001) and indicate this
estimate by $\b{\hat\eta}$;
\item[(ii)] convert $\b{\hat\eta}$ into $\b{\hat\pi}$ using the Newton-Rapson
algorithm described in Glonek \& McCullagh (1995) and Colombi \& Forcina
(2001);
\item[(iii)] choose
as importance density a Dirichlet distribution with mean vector equal to
$\b{\hat\pi}$, that is $g(\b\pi)\sim D(\alpha\b{\hat\pi})$, 
where $\alpha$ is a tuning
parameter that can be appropriately chosen so that enough draws from
the importance density satisfy the constraints imposed by model
$M_k$.  The optimal tuning parameter could be chosen by minimizing the
variance of the approximation, but this expression depends itself on
the target quantity. A simple approach, which we use in this paper, is
to try different values on a suitable grid (say from 0.02 to 50). 
\end{enumerate}
The same strategy can be adopted for choosing an appropriate
importance density to estimate $1/c_k$. In this case, the maximum
likelihood estimate $\b{\hat\eta}$ in (i) will be that 
corresponding to a hypothetical contingency table having a vector of
frequencies $\b y$ with all elements equal to zero.

To give an idea of the precision of the algorithm, we consider again
the above mentioned example in Section \ref{firstexample}. For those
data, we compare the true value of $1/c_k$, exactly computable for
the TP$_2$ model, with its estimates obtained in three separate runs
of the algorithm. The results, which are given in Table \ref{ckcheck},
show that the approximation is rather satisfactory in all
cases. 

\begin{table}[ht]\centering
\begin{tabular}{cccc}\hline\hline
True value & Estimate \#1 & Estimate \#2 &  Estimate \#3 \\
$2.9802\times 10^{-8}$ & $2.2315\times 10^{-8}$ & $2.8510\times 10^{-8}$
& $3.5776\times 10^{-8}$ \\
\hline\hline\end{tabular} \caption{\em True and estimated $1/c_k$
for model TP$_2$ on data in Section \ref{firstexample}.}
\label{ckcheck}
\end{table}

The second issue in estimating the Bayes factor arises in the
presence of about equality constraints. As already noted, for models
formulated accordingly to strict equality constraints, the Bayes
factor cannot be interpreted as the ratio between the proportions of
encompassing posterior and prior in agreement with the constraints,
since these proportions would be exactly zero. However, it has been
recently shown (Wetzels {\em et al.}, 2010) that the encompassing
approach naturally extends to exact equality constraints by
considering the ratio of the heights for the encompassing posterior
and prior distributions evaluated under the constraint (i.e., the
Savage-Dickey density ratio). However, this approach to handle
hypotheses specified through exact equality constraints complicates
the computation of the Bayes factor for models containing both
equality and inequality constraints. For this reason, we rather
preferred to follow the idea of Berger and Sellke (1987) and
Klugkist {\em et al.} (2010) of substituting exact equality with
about equality constraints. In this way, the interpretation of the
Bayes factor provided in (\ref{bf3}) is preserved and models
containing inequality or about equality constraints, as well as a
mix of both constraints, can be handled in a unified manner.
Moreover, Berger and Delampady (1987) noted that a Bayes factor
based on equality constraints is indistinguishable from a Bayes
factor based on about equality constraints, provided that the
interval around the exact equality constraint is small enough.
However, if this interval is too small, we incur again in the rare
event problem illustrated above, when trying to estimate $c_k^{-1}$
and $d_k^{-1}$.

To solve the above problem, Klugkist {\it et al.} (2010) proposed a
stepwise procedure which guarantees that a small enough interval is
used and does not actually need to pre-specify the size of this
interval. In principle, we could use the method of Klugkist {\it et
al.} (2010) to estimate the required constants $c_t$ and $d_t$. This
method is based on drawing random numbers from suitably truncated
Gamma distributions, which are then normalised to obtain the vector
$\b\pi$. The way in which the support of these distributions is
chosen depends on the adopted constraints. In our case, however, the
complexity of the constraints implies that it is difficult to define
how the support of each of these variables must be constrained; on
the other hand, a rejection sampling procedure to draw random values
from the truncated normal would be rather slow. For these reasons,
we prefer to adapt the iterative procedure in Klugkist {\it et al.}
(2010) exploiting, again, the importance sampling method.
According to our procedure, only two different samples, one drawn
from the importance density for the prior and the other one from the
importance density for the posterior, are required to estimate the
Bayes factor, thus overcoming the problem of sampling from
constrained distributions, which affects the procedure in Klugkist
{\em et al.} (2010). The details of the corresponding algorithm are
given in Appendix B.

Coming to the issue of parameter estimation, we need to acknowledge
that, for the class of models considered here, obtaining point or
interval estimates of the parameters is not in general of great
interest. The main interest rather lies in model selection as a tool
for evaluating which hypothesis is mostly supported by the data.
Nevertheless, once a particular model $M_k$ has been selected for
the data at hand, Bayesian parameter estimation is based on the
posterior distribution of the model parameters and, in our setting,
a sample from this posterior is already available after model choice
as a byproduct of the procedure to estimate $d_k$. In particular, we take
the set of all the draws from the parameter posterior distribution,
$D(\b 1_r+\b y)$, of the saturated model that are in agreement with
the constraints imposed by model $M_k$. This set should contain
enough draws to be also used for parameter estimation purposes.

Obviously, parameter estimates can be obtained in the way described
above for models defined by about equality constraints but not for
models defined by exact equality constraints. As an alternative, if
estimates under exact equality constraints are required, the
parameterisation in $\b\eta$ can be used, after choosing an appropriate
prior distribution on this parameter vector, for example a Gaussian
distribution. However, as already noticed, such a parameterisation would
complicate model selection in the presence of models expressed
through different types of logit.
\section{Applications}\label{exa}
In the following, we illustrate the proposed approach through three
applications involving some interesting datasets which also include
explanatory variables. In the first application, illustrated in
Section \ref{firstexample}, we also propose an analysis of
sensitivity with respect to the prior specification.
\subsection{Classification of men by social class and social class of their
fathers}
\label{firstexample}

We first consider a dataset (see Table \ref{father}) referred to a
sample of British males cross-classified according to their
occupational status ($A_2$) and that of their father ($A_1$).

\begin{table}[ht]\centering
\begin{tabular}{ccrrrrrr}\hline\hline
&& \multicolumn6c{$A_2$} \\
\cline{3-8} $A_1$
&&\multicolumn1c{I}&\multicolumn1c{II}&\multicolumn1c{
III}&\multicolumn1c{ IV}&\multicolumn1c{ V}&\multicolumn1c{
VI}\\\hline { I}&&{ 125}&{ 60}&{ 26}&{ 49}&{ 14}&{ 5}\\{ II}&&{
47}&{ 65}&{ 66}&{ 123}&{ 23}&{ 21}\\{ III}&&{ 31}&{ 58}&{ 110}&{
223}&{ 64}&{ 32}\\{ IV}&&{ 50}&{ 114}&{ 185}&{ 715}&{ 258}&{
189}\\{ V}&&{ 6}&{ 19}&{ 40}&{ 179}&{ 143}&{ 71}\\{ VI}&&{ 3}&{
14}&{ 32}&{ 141}&{ 91}&{ 106}\\\hline\hline\end{tabular}
\caption{\em Father ($A_1$) and son ($A_2$) occupational status
for a sample of 3,488 British males.} \label{father}
\end{table}

\newpage
The data have been already analysed by other authors, such as Goodman
(1991) and Dardanoni \& Forcina (1998). In particular, Dardanoni \& 
Forcina (1998), following a likelihood ratio approach, concluded that the data
conform to some forms of positive association. However, due to
presence of nuisance parameters, given by marginal column
probabilities, they did not reach a definitive conclusion about
TP$_2$.

For these data we first compared the saturated model ($M_1$) with
the independence model ($M_2$), the saturated model incorporating
PQD ($M_3$) and that incorporating TP$_2$ ($M_4$). For each of the
latter three models we estimated the Bayes factor with respect to
the saturated model, taken as reference model, through the 
algorithm described in Section \ref{modsel}. We obtained the
following results:
\begin{center}
\begin{tabular}{ccc}
 $\log(\hat{B}_{21})$ & $\log(\hat{B}_{31})$ & $\log(\hat{B}_{41})$\\
\hline -34.88 & 4.32 & 5.12
\end{tabular}
\end{center}

In order to give to the reader an idea of the computational details, we
point out that to compute $\hat B_{21}$ (see Section \ref{compbf}) we used
$\alpha=20$ for the prior and $\alpha=1$ for the posterior, where
as importance density we used a Dirichlet with parameters
corresponding to the independence model itself. Since we have about
equality constraints, we have used the algorithm in Appendix B starting
from $\epsilon=0.1$, with tuning parameter $b=0.5$. 
The algorithm has stopped after two iterations,
hence with $\epsilon=0.025$. The approximation has been replicated
$B=100$ times and we found it fairly stable (we obtained a
standard deviation of the 100 replicates smaller than 2). We report
the average estimate, which can be seen as a single estimate obtained
from a concatenated sample. The other two Bayes factors do not
involve about equality, but only inequality constraints. For the case of
PQD (i.e., $\hat B_{31}$) we did not need importance sampling because
sampling directly from the prior and posterior gives a large number of samples
satisfying the constraints. 

The hypothesis of independence must be definitely rejected, whereas
that of positive association may be accepted. In particular, the
model incorporating TP$_2$, formulated by requiring that all the
local log-odds ratios are non negative, has to be preferred to
that incorporating PQD, which is formulated through global log-odds
ratios. This means that the data conform to the strongest notion
of positive association among those considered by Douglas {\em et
al.} (1991). Hence, we can state that sons coming from a better
family have a higher chance of success also conditional on
remaining within any given subset of neighbouring classes. On the other hand,
the hypothesis of uniform association has to be rejected since,
comparing the model incorporating this constraint in addition to
TP$_2$ ($M_5$) with model $M_4$, we obtained
$\log(\hat{B}_{54})=-28.01$.

In order to perform some sensitivity analysis, we also calculated
the Bayes factors in the table above under other three different
Dirichlet prior parameters, obtaining the following results:
\begin{center}
\begin{tabular}{cccc}
& $\log(\hat{B}_{21})$ & $\log(\hat{B}_{31})$ & $\log(\hat{B}_{41})$\\
\hline
$D(0.5\b 1_r)$ & -34.49 & 4.26 & 5.04\\
$D(2\b 1_r)$ & -34.22 & 4.36 & 5.26\\
$D(5\b 1_r)$ & -32.18 & 4.39 & 6.04
\end{tabular}
\end{center}
It can be seen that there is only a slight sensitivity to prior
assumptions for these data. By varying the prior we do not reach
different conclusions with respect to model choice. We obtained
similar results, not reported here, for the other Bayes factors
computed in this section and in the next two sections.

Moving back to the data, we also considered some constraints on the
marginal distributions of the response variables. In particular, we
considered model $M_6$, formulated by incorporating in $M_4$ the
constraint that the marginal distributions of $A_1$ and $A_2$ are
equal, and model $M_7$, by incorporating in $M_4$ the constraint
that every local logit of $A_2$ is greater than the corresponding
local logit of $A_1$; this in turn implies that the marginal
distribution of $A_2$ is stochastically greater than that of $A_1$.
The Bayes factors of these two models with respect to $M_4$ are:
\begin{center}
\begin{tabular}{cc}
 $\log(\hat{B}_{64})$ & $\log(\hat{B}_{74})$\\ \hline  -0.78 & 2.22
\end{tabular}
\end{center}
Model $M_7$ seems to be supported by the data. This means that we
can observe not only {\em pure mobility}, that is positive
association between family's origin and the son's status, but also
{\em structural mobility}, which instead refers to how far apart
the two marginal distributions are and is essentially related to
socioeconomic growth.
\subsection{Classification of elderly people by Alzheimer's
disease and cognitive impairment}
The second dataset we analysed (see Table \ref{elderly}) is referred
to a sample of elderly people cross-classified by Alzheimer's
disease ($A_1$) and cognitive impairment ($A_2$), stratified by age
($B$); the data are taken from Agresti (1990, p. 298). The levels of
$A_1$ are: (IV) highly probable; (III) probable; (II) possible; (I)
unaffected; the levels of $A_2$ are: (V) severe; (IV) moderate;
(III) mild; (II) borderline; (I) unaffected.
\begin{table}[ht]
\centering
\begin{tabular}{ccrrrrrrrrrrr}\hline\hline
&&\multicolumn4c{$A_2$ ($< 75$)}&&\multicolumn4c{$A_2$ ($\geq 75$)}\\
\cline{3-6} \cline{8-11} $A_1$ && IV & III & II & I && IV & III &
II & I \\\hline { V}&&{ 2}&{ 1}&{ 1}&{ 0}&&{ 14}&{ 24}&{ 2}&{
0}\\{ IV}&&{ 1}&{ 12}&{ 10}&{ 1}&&{ 19}&{ 48}&{ 25}&{ 0}\\{
III}&&{ 0}&{ 8}&{ 27}&{ 5}&&{ 1}&{ 25}&{ 63}&{ 4}\\{ II}&&{ 0}&{
0}&{ 20}&{ 4}&&{ 0}&{ 0}&{ 35}&{ 7}\\{ I}&&{ 0}&{ 0}&{ 0}&{
85}&&{ 0}&{ 0}&{ 0}&{ 69}\\\hline\hline
\end{tabular} \caption{\em Alzheimer's
disease ($A_1$) and cognitive impairment ($A_2$) for a sample of
513 elderly people, stratified by age ($B$: less than 75, more
than 75).} \label{elderly}
\end{table}

As the categories of both response variables are in reverse order,
we based our analysis on reverse continuation logits. In this
setting, we compared the saturated model ($M_1$) with the model of
conditional independence ($M_2$) and the saturated model
incorporating positive association in every stratum ($M_3$). The
estimated Bayes factors are:
\begin{center}
\begin{tabular}{cc}
$\log(\hat{B}_{21})$ & $\log(\hat{B}_{31})$\\ \hline  -6.31 &
4.76
\end{tabular}
\end{center}
The hypothesis of conditional independence is not supported by the
data, whereas that of positive association in each stratum is
strongly supported. This means that, also conditionally on the age,
worst diagnoses of Alzheimer's disease are associated with most
severe cognitive impairment for both.

Then, we tried to test further hypotheses on the association between
the two response variables. In particular, we considered the following
constraints: (i) the level of the association is the same in each
stratum; (ii) the association is stronger in the first stratum;
(iii) the association is stronger in the second stratum. The models
obtained by incorporating these hypotheses in $M_3$ are denoted,
respectively, by $M_4$, $M_5$ and $M_6$. The estimated Bayes factors
for these models with respect to $M_3$ are:
\begin{center}
\begin{tabular}{ccc}
$\log(\hat{B}_{43})$ & $\log(\hat{B}_{53})$ &
$\log(\hat{B}_{63})$\\\hline -4.26 & -22.40 & 8.12
\end{tabular}
\end{center}
A certain amount of evidence in favour of model $M_6$ is noted.
Using this as reference model, further constraints on the
marginal distributions can be added, such as: the marginal
distribution of $A_1$ increases, namely the reverse continuation
logits decreases, with age ($M_7$); the marginal distribution of
$A_2$ increases with age ($M_8$); both marginal distributions
increase with age ($M_9$). We have the following results:
\begin{center}
\begin{tabular}{ccc}
$\log(\hat{B}_{76})$ & $\log(\hat{B}_{86})$ &
$\log(\hat{B}_{96})$\\ \hline 6.37 & 17.17 & 22.14
\end{tabular}
\end{center}
Models $M_7$, $M_8$ and $M_9$ seem to be all compatible with the
data, therefore we chose the one with the highest Bayes factor as
the most plausible one, that is model $M_9$. This implies that, as
age increases, individuals are more likely to have a serious level
of cognitive impairment and to be diagnosed the Alzheimer's disease
with a higher degree of confidence. Therefore, age does not only
affect the association between the two response variables, that is
stronger for elder people, but also shows a direct effect on their
marginal distributions.
\subsection{Clinical trial for skin disorder}
The dataset in Table \ref{trial}, already analyzed by Glonek \&
McCullagh (1995) and Koch {\it et al.} (1991), refers to a
clinical trial which, for confidentiality, was fictitiously
described as pertaining to the treatment of a skin disorder. The 72
subjects in the study are divided into two groups, the first one
receiving the treatment and the second one receiving the placebo. An
ordinal response variable, with levels poor/fair (I), good (II) and
excellent (III), was recorded for each subject on four different
occasions: 3 days ($A_1$), 7 days ($A_2$), 10 days ($A_3$) and 14
days ($A_2$) after treatment. Given the nature of the response
variables, it is natural to use global logits.
\begin{table}[ht]\centering
\begin{tabular}{cccccccccccccccccccccccccc}\hline\hline
& & & \multicolumn{9}{c}{Treatment}&&\multicolumn{9}{c}{Placebo}\\
\cline{4-12} \cline{14-22} & & $A_3$ & \multicolumn3c{I} &
\multicolumn3c{II} & \multicolumn3c{III} && \multicolumn3c{I} &
\multicolumn3c{II} &
\multicolumn3c{III}\\\cline{4-12}\cline{14-22}
 $A_1$ & $A_2$ & $A_4$ & I & II & III & I & II & III & I & II & III && I & II & III &
I & II & III & I & II & III\\ \hline
  & I && 0 &  1  & 0 & 0& 0 & 0 & 0 &  0 &  0 && 0  &   6  &
1 &  0  &   2 & 0 & 0& 0 &0 \\
I & II && 0  &   2  &   0 & 0  &   3  &   2  &   0 & 0 & 1 &&
0 &  3 & 1 & 0  & 6 & 2 & 0  &   0  &   0\\
  & III && 0 &   0  &   0  &   0   &  1  & 0   &  0 &    0  &   0 && 0   & 0
& 1 &  1 &   0  &   0  &   0 &  0   &  0 \\ & I && 0  & 0 & 0 &
0  &  3  &   0  & 0  &   0 &    0 && 0 &
0  & 0 &   0  & 0 &   0 &    0  &   0 &   0\\
II & II&& 0 & 0 & 1 & 0 & 2  & 4  &   1  &   1  & 0&& 0 & 1 & 0
& 0 & 1 & 2 & 0 & 3 & 3\\
& III && 0  & 0 & 0 &  0 & 1 & 3 & 0  & 0 & 5&& 0 & 0 &  0 & 0 &1
&0& 0  & 0 & 1 \\ &I&& 0 & 0 & 0 & 0 & 0 & 0 & 0 &0& 0&& 0 & 0
&   0 &  0 & 0 & 0  & 0 & 0 & 0\\
III&II&& 0 & 0  & 0 & 0 & 0  & 0 &   0  & 2  & 0&& 0 & 0 & 0
& 0 & 0 &   0 & 0 & 1 &  0 \\
&III&& 0 & 0 &  0  &  0  & 0   & 0   &  0   &  0 & 3&& 0 &    0&
0& 0  & 0 &   0  & 0 & 0 &  0
\\\hline\hline
\end{tabular} \caption{\em Response to treatment over time ($A_1$,
$A_2$, $A_3$, $A_4$) for a sample of 72 subjects, stratified by
type of treatment ($B$).} \label{trial}
\end{table}

The data are very sparse, as 128 of the 162 cells are
empty. Therefore, following Glonek \& McCullagh (1995), the largest
model we considered is a reduced model in which all the interactions
of order higher than two are set equal to zero.
Such a model, that we indicate by $M_2$,
is less restrictive than the largest model considered by Glonek \&
McCullagh (1995) which assumes that the association between every pair
of response variables is the same in the two strata.
Note that $M_1$, the saturated model, is still used as a reference
model to calculate the Bayes factors but is not included in the set
of models under choice.

We first compared model $M_2$ with the model that Glonek \&
McCullagh (1995) chose as final model ($M_3$). The latter is based
on the following constraints: (i) there is uniform association
within any stratum and between the strata; (ii) there is a constant
shift between marginal logits over time and between strata.
Comparing this model, with model $M_2$, we obtained
$\log(\hat{B}_{32})=0.19$. There is a very mild evidence in favor of
$M_3$. Finally, using $M_3$ as reference model, we considered model
$M_4$ obtained from $M_3$ by incorporating PQD and the constraint
that the marginal distribution of each response variable is
stochastically greater for the second stratum than for the first
one. These hypotheses seem to be supported by the data, as we have
$\log(\hat{B}_{43})=2.38$.
\section*{Acknowledgments}
F. Bartolucci and A. Farcomeni thank the `Einaudi Institute for Economics
and Finance' (EIEF) of Rome (IT) for financial support. 
\section*{References}
{\small \begin{description}

\item Agresti, A. (1990), {\it Categorical data Analysis},
John Wiley \& Sons, New York.

\item Albert, J. H. (1996), Bayesian selection of log-linear models,
{\em Canadian Journal of Statistics}, {\bf 24}, pp. 327-347.

\item Albert, J. H. (1997), Bayesian testing and estimation of
association in a two-way contingency table, {\em Journal of the
American Statistical Association}, {\bf 92}, pp. 685-693.

\item Bartolucci, F., Colombi, R. and Forcina, A. (2007), An extended class of marginal link
functions for modelling contingency tables by equality and inequality constraints, {\em
Statistica Sinica}, {\bf 17}, pp. 691-711.

\item Bartolucci, F., Forcina, A. \& Dardanoni, V. (2001), Positive
Quadrant Dependence and Marginal Modelling in Two-Way Tables with
Ordered Margins, {\it Journal of the American Statistical
Association}, {\bf 96}, pp. 1497-1505.

\item Berger, J.O. \& Delampady, M. (1987), Testing precise hypotheses,
  {\em Statistical Science}, {\bf 2}, pp. 317-352.

\item  Berger, J.O. \& Sellke, T. (1987), Testing a point null
hypothesis: The irreconcilability of p values and evidence, {\em
Journal of the American Statistical Association}, {\bf 82}, pp.
112-122.

\item Bergsma, W. \& Rudas, T. (2002), Marginal Models for categorical Data,
{\em The Annals of Statistics}, {\bf 30}, pp. 140-159.

\item Bucklew, J. (2004) {\em Introduction to rare event simulation},
  Springer series in Statistics, New York.



\item Colombi, R. \& Forcina, A. (2001), Marginal regression models for the
analysis of positive association of ordinal response variables,
{\it Biometrika}, {\bf 88}, pp. 1007-1019.


\item Dardanoni, V. \& Forcina, A. (1998), A Unified Approach to
Likelihood Inference on Stochastic Orderings in a Nonparametric Context,
{\it Journal of the American Statistical Association}, {\bf 93},
pp. 1112-1123.

\item Dellaportas, P. \& Forster, J. J. (1999), Markov chain Monte
Carlo Model Determination for Hierarchical and Graphical
Log-linear models, {\em Biometrika}, {\bf 86}, pp. 615-633.

\item Dawid, A.P. \& Lauritzen, S.L. (2001). Compatible prior
  distributions. {\em Bayesian Methods with Applications to Science, Policy
  and Official Statistics. Selected Papers from ISBA 2000: The Sixth
  World Meeting of the International Society for Bayesian Analysis}
(George, E.I. editor), Eurostat, Luxembourg, pp. 109-118.

\item Douglas, R., Fienberg, S. E., Lee, M. T., Sampson, A. R. \& Whitaker,
L. R. (1991), Positive Dependence Concepts for Ordinal Contingency
Tables, {\em Topics in Statistical Dependence} (H. W. Block, A. R.
Sampson \& T. H. Savits editors), Hayward, CA, pp. 189-202.

\item Glonek, G. (1996), A Class of Regression Models for Multivariate
Categorical Responses, {\em Biometrika}, {\bf 83}, pp. 15-28.

\item Glonek, G. \& McCullagh, P. (1995), Multivariate logistic
models, {\em Journal of the Royal Statistical Society, Series B},
{\bf 57}, pp. 533-46.

\item Goodman, L. (1991), Measures, Models and Graphical Displays
in the Analysis of Cross-Classified Data, {\it Journal of the
American Statistical Association}, {\bf 86}, pp. 1085-1111.

\item Green, P. J. (1995), Reversible jump Markov chain Monte
Carlo computation and Bayesian model determination, {\em
Biometrika}, {\bf 82}, pp. 711-732.

%
\item Hoijtink, H., Klugkist, I. \&
Boelen, P. (2008), {\em Bayesian Evaluation of Informative
Hypotheses}, Springer, New York.

\item Jeffreys, H. (1935), Some Tests of Significance, Treated by Theory
of Probability, {\em Proceeding of the Cambridge Philosophycal
Society}, {\bf 31}, pp. 203-222.
\item Jeffreys, H. (1961), {\em Theory of Probability, 3rd edition},
Oxford University press.
\item Karlin, S. (1968), {\em Total positivity}, Stanford University Press.

\item Kass, R. E. \& Raftery, A.E. (1995), Bayes factors,
{\em Journal of the American Statistical Association}, {\bf 90},
pp. 773-795.

\item Klugkist, I. \& Hoijtink, H. (2007), The Bayes factor for
  inequality and about equality constrained models, {\em Computational
    Statistics and Data Analysis}, {\bf 51}, pp. 6367-6379.

\item Klugkist, I., Kato, B. \& Hoijtink, H. (2005a), Bayesian model
  selection using encompassing priors, {\em Statistica Neederlandica},
  {\bf 59}, pp. 57-69.

\item Klugkist, I., Laudy, O. \&
Hoijtink, H. (2005b), Inequality constrained analysis of variance: a
Bayesian approach, {\em Psychological Methods}, {\bf 10}, pp.
477-493.

\item Klugkist, I., Laudy, O. \& Hoijtink, H. (2010), Bayesian
  evaluation of inequality and equality constrained hypotheses for
  contingency tables, {\em Psychologial methods}, {\bf 15}, pp. 281-299.

\item Koch, G., Singer, J., Stokes, M., Carr, G., Cohen, S. \&
Forthofer, R. (1991), Some aspects of weighted least squares
analysis for longitudinal categorical data, in {\em Statistical
Models for Longitudinal Studies of Health} (J.H. Dwyer, M.
Feinlib, P. Lippert \& H. Hoffmeister editors), Oxford University
press, pp. 215-258.

\item Lehmann, E. L. (1966), Some Concepts of Dependence,
{\em The Annals of Mathematical Statistics}, {\bf 37}, pp.
1137-1153.

\item McCullagh P. \& Nelder, J. A. (1989), {\em Generalized linear models, 2nd
edition}, Chapman and Hall, London.


\item Plackett, R. L. (1965), A Class of Bivariate Distributions,
{\em Journal of the American Statistical Association}, {\bf 60},
pp. 516-522.

\item Tuyl, F., Gerlach, R. \& Mengersen, K. (2009) Posterior
  predictive arguments in favor of the Bayes-Laplace prior as the
  consensus prior for binomial and multinomial parameters, {\em
    Bayesian Analysis}, {\bf 4}, pp. 151-158.

\item Wetzels, R., Grasman, R.P.P.P. \& Wagenmakers, E.-J. (2010) An
  encompassing prior generalization of the Savage-Dickey density
  ratio, {\em Computational Statistics \& Data Analysis}, {\bf 54},
pp. 2094-2102.

\end{description}}

\section*{Appendix}

\subsection*{A: Transformation from $\b\pi$ to $\b\eta$}
\label{matrices}
The matrices $\b C$ and $\b M$ in (\ref{marg}) may be obtained as
follows. $\b C$ is a block diagonal matrix with blocks $\b C_{\bl
z}$, ordered as $\b\eta_{\bl z}$ in $\b\eta$, given by
\[
\b C_{\bl z} = \bigotimes_{i=1}^q \b C_i,
\]
where $\b C_i = 1$ if $z_i=1$ and $\b C_i = \pmatrix{\b I_{m_i-1}
& -\b I_{m_i-1}}$ otherwise. Similarly, $M$ has blocks of columns
$M_z$ given by
\[
\b M_z = \bigotimes_{i=1}^q \b M_i,
\]
where $\b M_i = \b 1_{m_i}\tr$ if $z_i=0$; otherwise, we have
\[
\b M_i=\left\{\begin{array}{ll} \pmatrix{\b I_{m_i-1} & \b
0_{m_i-1} \cr \b 0_{m_i-1} & \b I_{m_i-1}} &  \mbox{ if logits of
type $l$ are used for the $i$-th variable,} \cr\vspace{-3mm}\cr
\pmatrix{\b T_{m_i-1} & \b 0_{m_i-1} \cr \b 0_{m_i-1} & \b
T_{m_i-1}\tr} & \mbox{ if logits of type $g$ are used for the
$i$-th variable,} \cr\vspace{-3mm}\cr \pmatrix{\b I_{m_i-1} & \b
0_{m_i-1} \cr \b 0_{m_i-1} & \b T_{m_i-1}\tr} & \mbox{ if logits
of type $c$ are used for the $i$-th variable,}\cr\vspace{-3mm}\cr
\pmatrix{\b T_{m_i-1} & \b 0_{m_i-1} \cr \b 0_{m_i-1} & \b
I_{m_i-1}} & \mbox{ if logits of type $r$ are used for the $i$-th
variable},
\end{array}\right.
\]
where $\b T_h$ is a $h\times h$ lower triangular matrix of ones.

\subsection*{B: Computing Bayes Factors with about equality constraints}
\label{bfeq}
First of all, we recall that about equality constraints are
specified as $|\b E\b \eta| \leq \b \epsilon$, for a small
$\b\epsilon>0$. If $\b\epsilon$ is too large, the corresponding
Bayes factor is far from the Bayes factor which would be obtained
with precise equality constraints. If $\b\epsilon$ is too small,
estimates of the proportion of the encompassing prior and
encompassing posterior in agreement with the constraints may be
inefficient.

In order to fix a suitable value for $\b\epsilon$, we adapt the
iterative procedure of Klugkist {\it et al.} (2010). Suppose we want
to estimate $B_{k1}$ for the constrained model $M_k$ versus the
encompassing model, where the constrained model is subject to $|\b
E\b \eta| \leq \b \epsilon$ and, possibly, $\b U\b \eta \geq \b 0$.
Our procedure comprises the following steps:
\begin{enumerate}
\item choose
a small value $\b \epsilon_1$ and define $M_{k.1}$ as the model $M_k$ in which 
$\b \epsilon$ is put equal to $\b \epsilon_1$;
\item estimate
$\hat B_{(k.1)1}=\hat c_{k.1}/\hat d_{k.1}$, where $\hat c_{k.1}^{-1}$ and
$\hat d_{k.1}^{-1}$ are, respectively, the proportions of the sample from the 
encompassing prior and posterior distributions in agreement with the constraints imposed by $M_{k.1}$;
\item define
$\b \epsilon_{2}=b\b \epsilon_{1}$, with $0 < b < 1$, and $M_{k.2}$ as the model
$M_k$ in which $\b \epsilon$ is put equal to $\b \epsilon_{2}$;
\item estimate
$\hat B_{(k.2)(k.1)}=(\hat c_{k.2}/\hat d_{k.2})/(\hat c_{k.1}/\hat d_{k.1})$,
where $\hat c_{k.2}^{-1}$ and $\hat d_{k.2}^{-1}$ are, respectively, the
proportions of the
samples from the encompassing prior and posterior in agreement with the 
constraints imposed by
$M_{k.2}$.
\end{enumerate}
Repeat steps 3 and 4, with each $\b \epsilon_{n+1}=b\b
\epsilon_{n}$, until the condition $\hat B_{(k.n+1)(k.n)}\approx 1$
is not satisfied. Then the required Bayes factor estimate $\hat
B_{k1}$ can be calculated by multiplication:
\begin{equation}\label{mult}\hat B_{k1} = \hat B_{(k.1)1} \times \hat
B_{(k.2)(k.1)} \times \cdots \times \hat B_{(k.n)(k.n-1)}.
\end{equation}
In the limit (i.e., when $\b \epsilon_n \to \b 0$), this method
yields the estimate of the Bayes factor for model $M_k$ with exact
equality constraints versus the encompassing model.

Notice that, in the procedure above, the problem of getting
inefficient estimates for the proportion of encompassing prior and
posterior in agreement with the constraints  is solved by using the
importance sampling approach described in Section \ref{compbf}.
Thus, only two different samples, one drawn from the importance
density for the prior and the other one from the importance density
for the posterior, are required to compute all the Bayes factor
estimates in (\ref{mult}).

\end{document}